\documentclass{birkjour}
\newtheorem{theorem}{Theorem}
\newtheorem{problem}[theorem]{Problem}
\newtheorem{definition}[theorem]{Definition}
\newtheorem{lemma}[theorem]{Lemma}
\newtheorem{conj}[theorem]{Conjecture}
\newtheorem{prop}[theorem]{Proposition}

\begin{document}

%-------------------------------------------------------------------------
% editorial commands: to be inserted by the editorial office
%
%\firstpage{1} \volume{228} \Copyrightyear{2004} \DOI{003-0001}
%\seriesextraline{This is the Concrete Title of this Book\br H.E. R and S.T.C. W, Eds.}
%\firstpage{1}

\title{A solution to two old problems by Menger concerning angle spaces}

%----------Author 1
\author{Luis Felipe Prieto-Mart\'inez}
%\author[Birkh\"auser]{Birkh\"{a}user Publishing Ltd.}

\address{%
Departamento de Matem\'aticas ETSI Industriales-INEI,\\ Universidad de Castilla-La Mancha, Ciudad Real,
Spain}

\email{luisfelipe.prieto@uclm.es}

%\thanks{This work was completed with the support of our

%----------Author 2
%\author{A Second Author}
%\address{The address of\br
%the second author\br
%sitting somewhere\br
%in the world}
%\email{dont@know.who.knows}



%\date{January 1, 2004}

\begin{abstract}
 Around 1930, K. Menger expressed his interest in the concept of abstract angle function. He introduced a general definition of this notion for metric and semi-metric spaces. He also proposed two problems concerning conformal embeddability of spaces endowed with an angle function into Euclidean spaces. These problems received attention in later years but only for some particular cases of metric spaces. In this article we first update the definition of angle function to apply to the larger class of spaces with a notion of betweenness, which seem to us a more natural framework.  In this new general setting, we solve the two problems proposed by Menger.
\end{abstract}

%%% ----------------------------------------------------------------------
\maketitle
%%% ----------------------------------------------------------------------
%\tableofcontents

\section{Introduction}
 In relation to his work about embeddability of metric spaces in $\mathbb E^n$,  Menger explored in  \cite{M.A} (among other articles) the  notion of \emph{angle} for a metric or semi-metric space $(\mathcal X,d)$. Let us denote by $\mathcal X^{3*}=\{(A,B,C)\in\mathcal X^3: B\neq A,C\}$.  He defined  \emph{angle functions} to be maps $\angle:X^{3*}\to \mathbb [0,\pi]$ such that 
\begin{equation} \label{eq.defiangle}\angle(A,B,C)=\angle(C,B,A)=\begin{cases}0 &\text{if }d(B,A)+d(A,C)=d(B,C) \text{ or } \\ &\;\;\;d(B,C)+d(C,A)=d(B,A)\\  \pi & \text{if }d(A,C)=d(A,B)+d(B,C) \\ >0 & \text{in other case} \end{cases}\end{equation}

\noindent For him, an \emph{angle space} was  a triple $(\mathcal X,d,\angle)$, where $X$ is a set, $d$ a distance or semi-distance and $\angle$ an \emph{angle function}. He also  stated the following two problems:

%Menger already realized that for his triples $(\mathcal X,d,\angle)$ it is interesting to discuss if the distance function and the angle function are compatible, in the sense that they are in the euclidean space, that is,  for every $(B,A,C)\in\mathcal X^{3*}$ we have:
%\begin{equation} \label{eq.ec} d(B,C)^2=d(B,C)^2+d(A,B)^2+2\cdot d(B,C)\cdot d(A,B)\cdot \cos\angle(B,A,C) \end{equation}

%\noindent In this article we use the term \emph{euclidean-compatible} for a distance and angle functions satisfying the condition above. \emph{Euclidean-compatiblity} is necessary for the triple $(\mathcal X,d,\angle)$ to be, both isometrically and conformally, embedded in the euclidean space $\mathbb E^n$.

\begin{problem}[literal copy from page 750 in \cite{M.A}] \label{question.mengerplane} It would be desirable to know relations between the angles of an angle-space, characteristic for those angle-spaces which are conformal to a subset of the plane and to the plane itself. 

\end{problem}

\begin{problem}[literal copy from page 750 in \cite{M.A}] \label{question.mengerec} In a metrical space (...) each three points are congruent to three points of the plane and, therefore, the definition of distances induces a definition of angles (viz. the angles in the corresponding plane triangles). A problem (connected, of course, with the first mentioned) is to find conditions which are necessary and sufficient in order that in an angle-space, which at the same time is a metrical space, the angles are identical with the angles induced by the definition of distance.

\end{problem}

Let us clarify that, for the author, in a metric space $(X,d)$ the angles induced by the definition of distance are those satisfying:
\begin{equation} \label{eq.ec}\forall A,B,C\in\mathcal X,\; d(B,C)^2=d(A,C)^2+d(A,B)^2+2\cdot d(A,C)\cdot d(A,B)\cdot \cos\angle(B,A,C) \end{equation}

After the work by Menger, other attempts have been done to provide a useful definition of \emph{angle space}, always, as far as we know, in the setting of metric spaces. On the one hand the more complicated the structure of the metric space is and the more axioms  are required for the angle function, the more desirable properties are satisfied by the angle space. So, some  years after the article by Menger,  W. A. Wilson in \cite{W1} studied the concept of \emph{angle} for convex complete metric spaces $\mathcal X$ satisfying what he called the ``four point property''. In this setting, there is a notion of line, ray and segment, and it is possible to define angles between them, not only for point triples. The authors of  \cite{VW} decided to restrict themselves to the study of  angles in linear normal spaces.  On the other hand, in the recent articles \cite{HM, M.AR} again the notion of angle between three points for general metric spaces is discussed, with the intention of providing a definition that coincides with the usual angle between geodesics if the space is a Riemannian manifold.

In this paper, we provide a simple and general axiomatic definition of \emph{angle function}  for spaces $\mathcal X$ that are not necessary metric or semi-metric (Definition \ref{defi.angle}). This eliminates the hassle of discussing the relation between the distance and angle functions a priori. We defend that the minimum structure  required for $\mathcal X$ in order to have a not artificial notion of angle is a \emph{relation of betweenness}. More about these relations is said in Section \ref{sect.b}. Once we have settled our framework, the rest of the article is devoted to solving the two problems proposed by Menger (Sections \ref{sect.prob1} and \ref{sect.prob2}).  As it often happens with known results in the bibliography concerning embeddability of metric spaces in the Euclidean space $\mathbb E^n$, the characterization obtained in this article  in response to Problem \ref{question.mengerec} may be impractical for spaces with an infinite number of points. But it may be of great interest for finite spaces, and in relation to some problems in Applied Mathematics. In Section \ref{sect.final} we give a glimpse of how this simple concept of angle function may be of interest in many branches of Mathematics, from  machine learning to the classic Inscribed Square Problem. To conclude this paper, a list of three open problems, that can be regarded as a continuation of this work, is provided.

%%%%%%%%%%%%%%%%%%%%%%%%%%%%%%%%%%%%%%%%%
\section{Relation of betweeness. A new definition of angle space} \label{sect.b}

There exists recent publications studying spaces with a notion of \emph{betweeness} (see \cite{A.B}). There are different options to define this concept. In \cite{H.BG} some axiomatic approaches  for geometric use are discussed. Among all the possibilities, we are interested in the, not very restrictive, following one:

\begin{definition}[pseumetric  betweenness] A \emph{relation of betweennes} in a set $\mathcal X$ is  a ternary relation $\mathcal B$. If $(A,B,C)\in \mathcal B$ we say that $B$ is \emph{between} $A$ and $C$. And if  either $(B,A,C)\in\mathcal B$, $(A,B,C)\in \mathcal B$ or $(A,C,B)\in\mathcal B$ then we say that $A,B,C$ are \emph{collinear}. We ask moreover that the relation $\mathcal B$ satisfies:

\begin{enumerate}

\item[(B1)] If $(A,B,C)\in \mathcal B$, then $A,B,C$ are distinct.

\item[(B2)] If $(A,B,C)\in\mathcal B$, then $(C,B,A)\in\mathcal B$.

\item[(B3)] If $(A,B,C)\in\mathcal B$, then $(B,A,C)\notin\mathcal B$.

\item[(B4)]  If $(P_1,P_2,P_3),(Q_1,Q_2,Q_3)\in\mathcal B$ and these two triples of points have two points in common, then any possible set of three different points in $\{P_1,P_2,P_3\}\cup\{Q_1,Q_2,Q_3\}$ is collinear.

\end{enumerate}

\end{definition}

In \cite{A.B}  the concept of \emph{pseudometric betweenness} is introduced and the choice of the name is explained. Their last axiom  is slightly different to ours. The reason of this change is explained below. Also in that article \emph{lines} and \emph{segments} are defined. Here we prefer not to introduce more notation than needed. Instead,  (B4) gives us a notion of collinearity also for sets of four or more points: we say that such a set is collinear if any of its subset of  three points is collinear.

Every semi-metric space has a natural relation $\mathcal B$ of betweenness. We say that  $B$  is between $A,C$ if and only if $d(A,B)+d(B,C)=d(A,C)$. So we can appreciate that the  notion of betweenness was already captured in the definition by Menger \eqref{eq.defiangle}. Thus, for us, spaces with a relation of betweenness are the more natural and general framework to define angles, since, in some sense, collinearity is inseparably linked to  angles.

\begin{definition} \label{defi.angle} Let $\mathcal X$ be a space with a notion of betweenness $\mathcal B$. $\angle:\mathcal X^{3*}\to \mathbb [0,\pi]$ is an \emph{angle function}  if:

\begin{enumerate}

\item[(i)]  it is symmetric in the first and third coordinates $\angle(A,B,C)=\angle(C,B,A)$,

\item[(ii)]  $\displaystyle{\angle(A,B,C)=\begin{cases}0 &\text{if }A=C,\; (B,A,C)\in\mathcal B \text{ or }(B,C,A)\in\mathcal B\\
 \pi &\text{if } (A,B,C)\in\mathcal B\\
>0 &  \text{in other case}\end{cases}},$

\item[(iii)] (Axiom of Collinearity) if  $(B,C,D)\in\mathcal B$,  then for all $A\neq B$, $\angle(A,B,C)=\angle(A,B,D)$.

\end{enumerate}

\end{definition}

As the reader may have noticed, we have added one more requirement with respect to \eqref{eq.defiangle}. This  \emph{Axiom of Collinearity}  seems to us to be basic. We think that a possible reason for Menger to discard it is that  he was mainly  interested in  angle functions defined  in metric spaces with a compatible distance, where it is implicitely required.

For us, \emph{angle spaces} are triples $(\mathcal X,\mathcal B,\angle)$ where $\mathcal X$ is a set, $\mathcal B$ is a \emph{relation of betweenness} and $\angle$ is an angle function. If such an angle space admits a distance function $d$ such that $\angle$ is induced by $d$ as explained in the definition (see Equation \eqref{eq.ec})  we say that $d$ and $\angle$ are \emph{Euclidean-compatible}. As usual, we say that a function between two angle spaces is \emph{conformal}  if it preserves angles and the betweenness relation and that it is a \emph{conformal equivalence} if, in addition, it is biyective.

%%%%%%%%%%%%%%%%%%%%%%%%%%%%%%%%%%%%%%%%%%%%%%%%%%%
\section{Solution to Problem \ref{question.mengerec}} \label{sect.prob1}

An angle space $\mathcal X$ is called \emph{trivial} if it has two points or every three points in $\mathcal X$ are collinear. Trivial angle spaces obviously admit infinitely many Euclidean-compatible distance functions.  So we focus our study on the rest of the cases. In the following, for a given space $\mathcal X$ we use the terms \emph{trigon}, \emph{tetragon}, \emph{pentagon} and \emph{hexagon} to refer to subsets of  of $3, 4, 5$ and 6 different points in $\mathcal X$, respectively.

Let us begin with the following:

\begin{lemma} \label{lemm.nec} Let $(\mathcal X,\mathcal B,\angle)$ be an angle space. If this angle space admits a Euclidean-compatible distance function $d$, then for every trigon $\{A,B,C\}$  in $\mathcal X$:

\begin{enumerate}

\item[(EC1)] The sum of the angles equals $\pi$, that is:
$$ \angle(B,A,C)+\angle(A,B,C)+\angle(A,C,B)=\pi.$$

\end{enumerate}

\noindent Moreover, the angle space admits a Euclidean-compatible distance function $d$ if and only if for every non-trivial trigon $\{A,B,C\}$  in $\mathcal X$  the \emph{Law of Sines} holds: 
$$\frac{d(B,C)}{\sin\angle(B,A,C)}=\frac{d(A,C)}{\sin\angle(A,B,C)}=\frac{d(A,B)}{\sin\angle(A,C,B)}$$

\end{lemma}

\begin{proof} If $(\mathcal X,\mathcal B,\angle)$ admits a Euclidean-compatible distance function $d$ then for every trigon $\{A,B,C\}\subset\mathcal X$ the following system needs to have a solution in the indeterminates $a,b,c$.
\begin{equation}\label{eq.syscomp}\begin{cases}a^2=b^2+c^2-2bc\cos\alpha\\ b^2=a^2+c^2-2ac\cos\beta\\ c^2=a^2+b^2-2ab\cos\gamma\end{cases}\end{equation}

\noindent where $a=d(B,C)$, $b=d(A,C)$, $c=d(A,B)$, $ \alpha=\angle(C,A,B)$, $\beta=\angle(A,B,C)$ and $\gamma=\angle(B,C,A)$. In the following, let us assume that $\gamma\geq\alpha,\beta$.

 From the last equation we obtain, since $c>0$, that  $c=\sqrt{a^2+b^2-2ab\cos\gamma}$. Let us substitute this value in the other two equations:
$$\begin{cases}b^2-ab\cos\gamma=b\cos\alpha\sqrt{a^2+b^2-2ab\cos\gamma}\\ 
a^2-ab\cos\gamma=a\cos\beta \sqrt{a^2+b^2-2ab\cos\gamma}\\
c=\sqrt{a^2+b^2-2ab\cos\gamma}\end{cases} $$

\noindent Since $A,B,C$ are different, $a,b\neq 0$. Also, since the points are not collinear, $\sin\alpha,\sin\beta\neq 0$. So we can divide the first two equations by $ab$ and the third one by $a$ obtaining:
\begin{equation} \label{eq.magdalenafinal} \begin{cases}\frac{b}{a}-\cos\gamma=\cos\alpha\sqrt{(\frac{b}{a}-\cos\gamma)^2+\sin^2\gamma}\\ 
\frac{a}{b}-\cos\gamma=\cos\beta\sqrt{(\frac{a}{b}-\cos\gamma)^2+\sin^2\gamma}\\
\frac{c}{a}=\sqrt{(\frac{b}{a}-\cos\gamma)^2+\sin^2\gamma} \end{cases}\end{equation}

\noindent From the first two equations, we obtain respectively that $\frac{b}{a}-\cos\gamma=\sin^2\gamma\cot^2\alpha$, $\frac{a}{b}-\cos\gamma=\sin^2\gamma\cot^2\beta$, and then get:
$$\begin{cases}\frac{b}{a}=\frac{\sin(\alpha\pm\gamma)}{\sin\alpha}\\ \frac{a}{b}=\frac{\sin(\beta\pm\gamma)}{\sin\beta} \\ \frac{c}{a}=\frac{\sin\gamma}{\sin\alpha} \end{cases}$$

 Note that, in each of the two first equations, we have introduced an incorrect solution. Since $\gamma$ is the biggest angle and all of them are in $[0,\pi]$, the sign must be ``$+$'' in the first two  equations in the previous system. So we have already shown that the system in Equation \eqref{eq.syscomp} has a solution if and only if the following system does:
\begin{equation} \label{eq.revision2}\begin{cases}\frac{b}{a}=\frac{\sin(\alpha+\gamma)}{\sin\alpha}\\ \frac{a}{b}=\frac{\sin(\beta+\gamma)}{\sin\beta} \\ \frac{c}{a}=\frac{\sin\gamma}{\sin\alpha} \end{cases}\end{equation}

To prove (EC1), see that the system above has a solution if and only if $\sin(\alpha+\gamma)\sin(\beta+\gamma)=\sin \alpha \sin\beta$ and then $\cos(\alpha+\beta+ 2\gamma)=\cos(\alpha+\beta)$. Taking into account that $\alpha,\beta,\gamma\in[0,\pi]$,  this implies $ \pi=\alpha+\beta+\gamma $.

Now the Law of Sines being equivalent to the Euclidean-compatibility condition is a direct consequence of the fact that the system in Equation \eqref{eq.revision2} is, in turn, equivalent to the one in Equation \eqref{eq.syscomp}.

\end{proof}

 We are finally ready to prove the following theorem, which is a solution for Problem \ref{question.mengerec}:

%%%%%%%%%%%%%%%%%%%%%%%%%%%%%%%%%%%%%%%%%%%%%%%%%
%%%%%%%%%%%%%%%%%%%%%%%%%%%%%%%%%%%%%%%%%%%%%%%%
\begin{theorem} \label{theo.ec} Let $(\mathcal X,\mathcal B,\angle)$ be a non-trivial angle space with four or more points. Then, it admits a Euclidean-compatible distance function $d$ if and only if every subspace of four points does. This Euclidean-compatible distance function $d$ is  unique up to multiplication by a constant.

\end{theorem}

\begin{proof}  Let $(\mathcal X,\mathcal B,\angle)$ be a non-trivial angle space with four or more points and such that  every subspace of four points admits a Euclidean-compatible distance. Let $A,B\in\mathcal X$, and let $\lambda\in\mathbb R_{\geq 0}$. Let us construct  a function $d_{A,B,\lambda}:\mathcal X^2\to\mathbb R_{\geq 0}$ such that it is a distance, it is Euclidean-compatible to the angle space structure and such that $d_{A,B,\lambda}(A,B)=\lambda$.

First, let us prove that this distance function  $d_{A,B,\lambda}$, if it exists, is unique. This is an immediate consequence of the Law of Sines. In the following, we explain how this function is determined in each case:

%First, see that the fact that every trigon satisfies the Law of Sines, that $d_{A,B,\lambda}(A,B)=\lambda$  and that for every $P,Q\in\mathcal X$, $d_{A,B,\lambda}(P,P)=0$, $d_{A,B,\lambda}(P,Q)=d_{A,B,\lambda}(Q,P)$ and  $d_{A,B,\lambda}(P,Q)$ is linear in $\lambda$, already determine a unique function $d_{A,B,\lambda}$. 
%To do so, we need to check that it is possible to solve the problems in Figure 1 from this information.
%\begin{figure}[h]
%\includegraphics[width=\textwidth]{definiciond.pdf}
%\label{hola}
%\caption{In these diagrams, if three or more points lie in the same line, then they are collinear.}
%\end{figure}

\begin{itemize}

\item \emph{Case 1}. If $C$ is not collinear to $A,B$, we can use the Law of Sines in the trigon $\{A,B,C\}$ to obtain that:
\begin{equation} \label{eq.case1} d_{A,B,\lambda}(A,C)=\frac{\sin\angle(A,B,C)}{\sin\angle(A,C,B)} \lambda. \end{equation}

\noindent With a similar argument we can also compute $d_{A,B,\lambda}(B,C)$.

\item \emph{Case 2}. If  $C$ is collinear to $A,B$, let us choose a point  $E$ not collinear to $A,B,C$. Then, we can compute $d_{A,B,\lambda}(A,E)$ (using the Law of Sines  in the trigon $\{A,B,E\}$) and finally obtain  (using the Law of Sinces in the trigon $\{A,C,E\}$) that:
\begin{equation}\label{eq.case2} d_{A,B,\lambda}(A,C)= \frac{\sin\angle(A,E,C)}{\sin\angle(A,C,E)}\frac{\sin(A,B,E)}{\sin(A,E,B)}\lambda.\end{equation}

\noindent With a similar argument we can also compute $d_{A,B,\lambda}(B,C)$.

\item \emph{Case 3}.  Suppose that $D$ is not collinear to $A,B$. Use the Law of Sines  in $\{A,B,D\}$ and after this in $\{A,C,D\}$, we obtain that 
\begin{equation}\label{eq.pruebainfierno} d_{A,B,\lambda}(C,D)=\frac{\sin\angle(C,A,D)}{\sin\angle(A,C,D)}\frac{\sin\angle(A,B,D)}{\sin\angle(A,D,B)}\lambda.\end{equation}

\noindent If $D$ is collinear to $A,B$ but $C$ is not, a similar formula is derived exchanging the role of $C$ and $D$.

\item \emph{Case 4}. If $A,B,C,D$ are collinear, choose an element $E$ neither  collinear to $\{A,B\}$, nor to $\{C,D\}$. Then using again the Law of Sines (in the trigon $\{A,B,E\}$) we can compute  $d_{A,B,\lambda}(A,E)$  and then compute $d_{A,B,\lambda}(C,D)$ from this value as in the previous case:
\begin{equation}\label{eq.pruebainfiernocolineal1} d_{A,B,\lambda}(C,D)=\frac{\sin\angle(C,A,D)}{\sin\angle(A,C,D)}\frac{\sin\angle(A,E,D)}{\sin\angle(A,D,E)}\frac{\sin\angle(A,B,E)}{\sin\angle(A,E,B)}\lambda. \end{equation}

% $d_{A,B,\lambda}(A,E)=\frac{\sin(A,B,E)}{\sin(A,E,B)}\lambda$
\end{itemize}

Now, let us complete the argument for existence of such  function $d_{A,B,\lambda}$.

%%%%%%%%%%%%%%%%%%%%%%%%%%%%%%%%%%%%%%%%%%%%%%%%%%%%%%%
%\noindent \textbf{Claim 1.} \emph{In the definition above, for every $P,Q\in\mathcal X$:}
%$$d_{A,B,\lambda}(P,Q)=d_{B,A,\lambda}(P,Q)=d_{A,B,\lambda}(Q,P) $$
%The statement is obvious when we are in Cases 1 and 2. In Case 3, the statement is guaranteed by (EC2). According to Equations \eqref{eq.pruebainfierno} and \eqref{eq.susto}, the value in \eqref{eq.pruebainfierno} remains invariant if we exchange $A,B$ and $C,D$. Something similar happens to Case 4.

%%%%%%%%%%%%%%%%%%%%%%%%%%%%%%%%%%%%%%%%%%%%%%%%%%%%%%%
\noindent \textbf{Claim 1.} \emph{Let $\{A,B,C,D,E\}$ be a pentagon in $\mathcal X$. In this pentagon, there exists a distance function $d_{A,B,\lambda}$ which is Euclidean-compatible to the angle space structure and such that $d_{A,B,\lambda}(A,B)=\lambda$.}

\medskip

This statement holds automatically if all the points in the pentagon are collinear. Suppose that $E$ is not collinear to $A,B$.

Using the hypothesis in the theorem, the tetragon $\{A,B,C,D\}$ admits a Euclidean-compatible distance. Taking into account that if a given distance is Euclidean-compatible so are its multiples by a constant, we can see that  for every choice of $\lambda\in\mathbb R_{>0}$, there is a  Euclidean-compatible distance, denoted by $d_{A,B,\lambda}$, satisfying $d_{A,B,\lambda}(A,B)=\lambda$.

We only have to check that the distances between  $E$ and any of the points in $\{A,B,C,D\}$ are well-defined.

Using, again, the hypothesis in the theorem, the distance function is well-defined on $\{A,B,C,E\}$ and $\{A,B,D,E\}$. This gives a consistent definition of $d_{A,B,\lambda}(A,E)$, $d_{A,B,\lambda}(B,E)$, and $d_{A,B,\lambda}(C,E)$ in $\{A,B,C,E\}$ and of $d_{A,B,\lambda}(A,E)$, $d_{A,B,\lambda}(B,E)$ and $d_{A,B,\lambda}(D,E)$  in $\{A,B,D,E\}$.

The only thing that could spoil the coherence of the definition of the function $d_{A,B,\lambda}$ is that the values of $d_{A,B,\lambda}(A,E)$ and $d_{A,B,\lambda}(B,E)$ obtained for the tetragons $\{A,B,C,E\}$ and $\{A,B,D,E\}$ were different. But this cannot happen, since  both can be inferred from the non-trivial trigon $\{A,B,E\}$ which is contained in both tetragons.

%%%%%%%%%%%%%%%%%%%%%%%%%%%%%%%%%%%%%%%%%%%%%%%%%%%%%%%%%%%
\noindent \textbf{Claim 2.} \emph{Let $\{A,B,C,D,E,F\}$ be an hexagon in $\mathcal X$. In this hexagon, there exists a distance function $d_{A,B,\lambda}$ which is Euclidean-compatible to the angle space structure and such that $d_{A,B,\lambda}(A,B)=\lambda$.}

Again the statement is immediate if the hexagon is trivial, so we may assume that $F$ is not collinear to $A,B$. Using the previous claim, the pentagon $\{A,B,C,D,E\}$ admits a Euclidean-compatible distance function such that $d_{A,B,\lambda}(A,B)=\lambda$.

So we only need to check that the distances from $F$ to any other of the points in the pentagon $\{A,B,C,D,E\}$ are well-defined. Each of those distances appear in, at most, two of the pentagons containing $A,B$: $\{A,B,C,D,F\}$, $\{A,B,C,E,F\}$, $\{A,B,D,E,F\}$.

As in the previous case, the distances $d_{A,B,\lambda}(A,F)$ and $d_{A,B,\lambda}(B,F)$ are obtained from the trigon $\{A,B,F\}$, which belongs to the three pentagon listed above, and so they are well-defined.

To see that $d_{A,B,\lambda}(C,F)$ is well-defined, we need to check that the value obtained from the pentagons $\{A,B,C,D,F\}$ and $\{A,B,C,E,F\}$ coincide. To see this, note that it is determined by the tetragon $\{A,B,C,F\}$ which is common to both of them. A similar argument shows that $d_{A,B,\lambda}(D,F)$ is well-defined too.

%%%%%%%%%%%%%%%%%%%%%%%%%%%%%%%%%%%%%%%%%%%%%%%%%%%

\noindent \textbf{Claim 3.} \emph{For every $\lambda>0$, there exists a Euclidean-compatible distance function defined in $\mathcal X$ and such that $d_{A,B,\lambda}=\lambda$.}

First, we need to check that the function $d_{A,B,\lambda}$, as defined in Equations \eqref{eq.case1}, \eqref{eq.case2},\eqref{eq.pruebainfierno} and \eqref{eq.pruebainfiernocolineal1}, is well-defined.

The definition, does not depend on the choice of $E$ for Case 2. This is a consequence of the fact that the distance is well-defined in every pentagon containing $A,B$ (Claim 1). The definition, does not depend on the choice of $E$ for Case 4, either. This is a consequence of the fact that the distance  is well-defined in every hexagon containing $A,B$ (Claim 2).

Finally, we only need to prove that the function $d_{A,B,\lambda}$ is a Euclidean-compatible distance. The \emph{Identity of Indiscernibles} is easy to check from the definition of  $d_{A,B,\lambda}$. To prove the \emph{symmetry} in the two variables, we need to check that the choice of  trigons taken in Cases 2, 3 and 4 does  not affect the result. This is ensured by the fact that $d_{A,B,\lambda}$ is a well-defined distance function for every tetragon and pentagon in $\mathcal X$. This function $d_{A,B,\lambda}$ satisfies the condition of Euclidean-compatibility. We have this guaranteed by construction, since $d_{A,B,\lambda}$ is a Euclidean-compatible distance for every tetragon (and thus for every trigon) in $\mathcal X$. This also implies the \emph{Triangle Inequality}.

\end{proof}

%Finally, let us note that the distance function $d_{A,B,\lambda}$  is compatible with the betweenness relation $\mathcal B$ (thanks to (EC2)).

Finally, for completeness, we include the following characterization of non-trivial angle spaces  with three or four points admitting a Euclidean-compatible distance:

\begin{prop} Let $(\mathcal X,\mathcal B, \angle)$ be a non-trivial angle space.

\begin{itemize}

\item[(a)] If $\mathcal X$ consists of three points $A,B,C$, then it admits a Euclidean-compatible distance function if and only if it satisfies (EC1).

\item[(b)] If $\mathcal X$ consists of 4 points, such that three of them are collinear and the other one is not collinear to the rest, then it admits a Euclidean-compatible distance function if and only if every trigon $\{A,B,C\}\subset\mathcal X$ satisfies (EC1) and

\begin{itemize}

\item[(EC2)]   If  $(X_1,X_2,X_3)\in\mathcal B$ and $X_4$ is not collinear to $\{X_1,X_2,X_3\}$ then:
$$\begin{cases} \angle(X_1,X_2,X_4)+\angle(X_3,X_2,X_4)=\pi\\  \angle(X_1,X_4,X_2)+\angle(X_2,X_4,X_3)=\angle(X_1,X_4,X_3)\end{cases}.$$

\end{itemize}

\item[(c)] If  $\mathcal X$ does not contain collinear points  it admits a Euclidean-compatible distance function if and only if every trigon $\{A,B,C\}\subset\mathcal X$ satisfies (EC1) and

\begin{itemize}

\item[(EC3)] The following identity holds, as well as any other obtained from a permutation of the subindices:
$$\displaystyle{\sin\angle(X_1,X_2,X_3)\cdot \sin\angle(X_1,X_3,X_4)\cdot \sin\angle(X_1,X_4,X_2)=}$$
$$=\sin\angle(X_1,X_2,X_4)\cdot\sin\angle(X_1,X_4,X_3)\cdot\sin\angle(X_1,X_3,X_2) .$$

\item[(EC4)] The following identity holds, as well as any other obtained from a permutation of the subindices:
$$\frac{\sin\angle(X_1,X_2,X_3)}{\sin\angle(X_1,X_3,X_2)}\frac{\sin\angle(X_3,X_1,X_4)}{\sin\angle(X_3,X_4,X_1)}= $$
$$=\frac{\sin\angle(X_1,X_2,X_4)}{\sin\angle(X_1,X_4,X_2)}\frac{\sin\angle(X_4,X_1,X_3)}{\sin\angle(X_4,X_3,X_1)}=$$
$$=\frac{\sin\angle(X_2,X_1,X_3)}{\sin\angle(X_2,X_3,X_1)}\frac{\sin\angle(X_3,X_2,X_4)}{\sin\angle(X_3,X_4,X_2)}=$$
$$=\frac{\sin\angle(X_2,X_1,X_4)}{\sin\angle(X_2,X_4,X_1)}\frac{\sin\angle(X_4,X_2,X_3)}{\sin\angle(X_4,X_3,X_2)}.$$

\end{itemize}

\end{itemize}

\end{prop}

%%%%%%%%%%%%%%%%%%%%%%%%%%%%%%%%%%%%%%%%%%%%%%%%%%%%%%%
\begin{proof}  Part (a) is a consequence of Lemma \ref{lemm.nec}.

To prove (b), if such a distance function $d$ exists, let us use the following notation for the distance between points in $\mathcal X$:
$$x=d(X_1,X_2),\;x'=d(X_2,X_3),\;a=d(X_1,X_4),\;a'=d(X_3,X_4),\;h=d(X_2,X_4)$$

\noindent and let us denote the angles as:
$$\alpha=\angle(X_1,X_2,X_4),\;\alpha'=\angle(X_3,X_2,X_4),\;\beta=\angle(X_2,X_1,X_4)=\angle(X_3,X_1,X_4),$$
$$\beta'=\angle(X_2,X_3,X_4)=\angle(X_1,X_3,X_4),\;\gamma=\angle(X_1,X_4,X_2),\;\gamma'=\angle(X_2,X_4,X_3),$$
$$\delta=\angle(X_1,X_4,X_3).$$ 

\noindent Note that $d(X_1,X_3)=x+x'$, since $(X_1,X_2,X_3)\in\mathcal B$. If we use the Law of Sines simultaneously to the trigons $\{X_1,X_2,X_4\}$, $\{X_1,X_3,X_4\}$ and $\{X_2,X_3,X_4\}$, we get:
\begin{equation} \label{eq.syssac}\begin{cases}\frac{a}{\sin\alpha}=\frac{h}{\sin\beta}=\frac{x}{\sin\gamma} \\ \frac{a'}{\sin\alpha'}=\frac{h}{\sin\beta'}=\frac{x'}{\sin\gamma'}\\ \frac{x+x'}{\sin\delta}=\frac{a}{\sin\beta'}=\frac{a'}{\sin\beta}\end{cases}.\end{equation}

Let us prove the necessity of Condition (EC2). From the system above, we get the following one:
$$\begin{cases}\frac{\sin\gamma}{\sin \alpha}a+\frac{\sin\gamma'}{\sin\alpha'}a'=\frac{\sin\delta}{\sin\beta}a'\\ 
\frac{\sin\gamma}{\sin\alpha}a+\frac{\sin\gamma'}{\sin\alpha'}a'=\frac{\sin\delta}{\sin\beta'}a\end{cases}. $$

\noindent which should have a non-zero solution. So the determinant of the matrix of coefficients (once the equations are suitably fixed) must be zero. This yields:
$$\frac{\sin\gamma'}{\sin\alpha'\sin\beta'}+\frac{\sin\gamma}{\sin\alpha \sin\beta} =\frac{\sin\delta}{\sin\beta\sin\beta'}$$

\noindent and then, since $\sin\gamma=\sin(\alpha+\beta)$, $\sin\gamma'=\sin(\alpha'+\beta')$ we get:
$$\sin\beta\sin\beta'(\cot\alpha+\cot\alpha')+\sin(\beta+\beta')=\sin\delta $$

\noindent As $\delta+\beta+\beta'=\pi$, we have that $\sin(\beta+\beta')=\sin\delta$ and so $\cot\alpha=-\cot\alpha'$ and $\alpha=\pi-\alpha'$. Now since:
$$\begin{cases}\alpha+\beta+\gamma=\pi\\ \alpha'+\beta'+\gamma'=\pi\\ \beta+\beta'+\delta=\pi \end{cases}  $$

\noindent we can conclude that $\delta=\beta+\beta'$ (which is Condition (EC2)).

Finally, assuming that (EC2) holds, we have that \eqref{eq.syssac} is equivalent to:
$$\begin{cases}\frac{a}{\sin\alpha}=\frac{h}{\sin\beta}=\frac{x}{\sin\gamma} \\ 
\frac{a'}{\sin\alpha'}=\frac{h}{\sin\beta'}=\frac{x'}{\sin\gamma'}\\
\frac{a}{\sin\beta'}=\frac{a'}{\sin\beta}\end{cases}
$$

\noindent which has a unique solution up to multiples in the variables $a,a',x,x',h$.

To prove (c), let us construct a Euclidean-compatible distance $d$. Set $d(X_1,X_2)=\lambda$, for some $\lambda\in\mathbb R_{>0}$. Use the Law of Sines in the trigons $\{X_1,X_2,X_3\}$, $\{X_1,X_2,X_4\}$ to determine the values:
$$d(X_1,X_3)=\frac{\sin\angle(X_1,X_2,X_3)}{\sin\angle(X_1,X_3,X_2)}\lambda,\qquad d(X_2,X_3)=\frac{\sin\angle(X_2,X_1,X_3)}{\sin\angle(X_1,X_3,X_2)}\lambda$$
$$d(X_1,X_4)=\frac{\sin\angle(X_1,X_2,X_4)}{\sin\angle(X_1,X_4,X_2)}\lambda,\qquad d(X_2,X_4)=\frac{\sin\angle(X_2,X_1,X_4)}{\sin\angle(X_1,X_4,X_2)}\lambda$$

We need to check that the trigons $\{X_1,X_3,X_4\}$, $\{X_2,X_3,X_4\}$ satisfy the Law of Sines. We have that:
$$\frac{d(X_1,X_3)}{\sin\angle(X_1,X_4,X_3)}=\frac{d(X_1,X_4)}{\sin\angle(X_1,X_3,X_4)} $$
$$\frac{d(X_2,X_3)}{\sin\angle(X_2,X_4,X_3)}=\frac{d(X_2,X_4)}{\sin\angle(X_2,X_3,X_4)} $$

\noindent if and only if the equation in (EC3) holds (and one of its permutations).

From the trigons $\{X_1,X_3,X_4\}$, $\{X_2,X_3,X_4\}$ and the Law of Sines, we can obtain the value of $d(X_3,X_4)$ in different ways:
$$d(X_3,X_4)=\frac{\sin\angle(X_3,X_1,X_4)}{\sin\angle(X_1,X_4,X_3)}d(X_1,X_3)=\frac{\sin\angle(X_3,X_1,X_4)}{\sin\angle(X_1,X_3,X_4)}d(X_1,X_4)=$$
$$=\frac{\sin\angle(X_3,X_2,X_4)}{\sin\angle(X_2,X_4,X_3)}d(X_2,X_3)=\frac{\sin\angle(X_3,X_2,X_4)}{\sin\angle(X_2,X_3,X_4)}d(X_2,X_4).$$

\noindent All those values coincide if and only if the equation in (EC4) holds.

The role played by $X_1,X_2$ has been chosen arbitrarily, so the identities obtained permuting the subindices are also necessary.

\end{proof}

In the study of the trigonometry of the tetrahedron, the identities appearing in (EC3) are sometimes referred as \emph{Alternating Sines Theorem} and the identities appearing in (EC4) already appeared in the classic compendium of results  by G. Richardson \cite{R}.

%%%%%%%%%%%%%%%%%%%%%%%%%%%%%%%%%%%%%%%%%%%%%%%%%%%%%
%%%%%%%%%%%%%%%%%%%%%%%%%%%%%%%%%%%%%%%%%%%%%%%%%%%%%

%%%%%%%%%%%%%%%%%%%%%%%%%%%%%%%%%%%%%%%%%%%%%%%%%%%%%%%%
\section{Solution to Problem \ref{question.mengerplane}} \label{sect.prob2}

Menger  in \cite{M.F} and some years later C.  Morgan in \cite{M.E} provided characterizations of the metric spaces $(\mathcal X,d)$ that can be embedded in $\mathbb E^n$. Define for every $A,B,C\in \mathcal X$ the quantity $\langle A,B,C\rangle=\frac{1}{2}\left(d(A,C)^2+d(B,C)^2-d(A,B)^2\right)$. For each  $n$-simplex $\{X_0,\ldots,X_n\}\subset \mathcal X$ define 
$$\text{Vol}_n(\{X_0,\ldots,X_n\})=\frac{1}{n!}\sqrt{D_n(X_0,\ldots,X_n)} $$
where $D_n(X_0,\ldots,X_n)$ is the determinant of the  $n\times n$ matrix whose $(i,j)$ entry equals $\langle X_i,X_j,X_0\rangle$. As explained in \cite{M.E}, under some hypothesis this value is, in fact, the volume of the simplex. We say that a space $(\mathcal X,d)$ is $n$-\emph{flat} if for every $n$-simplex $\text{Vol}_n(\{X_0,\ldots,X_n\})$ is a real number. In this case, its  \emph{dimension}  is the largest integer (if it exists) for which there is an $n$-simplex for which $\text{Vol}_n(\{X_0,\ldots,X_n\})>0$.  We have that:

\begin{theorem}[proved by Morgan in \cite{M.E}]  \label{theo.morgan}

A metric space $(\mathcal X,d)$ can be embedded in $\mathbb E^n$ if and only if it is $k$-flat for every $1\leq k\leq n$, and of dimension less than or equal to $n$.

\end{theorem}

 In relation to our setting, if the metric space $(\mathcal X,d)$ has also an angle function that is Euclidean compatible with $d$, then $\langle A,B,C\rangle=d(A,C)\cdot d(B,C)\cdot \cos(\angle(A,C,B))$ and  so
\begin{equation} \label{eq.cosenoso} D_n(X_0,\ldots,X_n)=(\Pi_{j=1}^n d(X_0,X_j)^2)\begin{array}{| c c c |} \cos\angle (X_1,X_0,X_1) & \hdots & \cos \angle (X_1,X_0,X_n)\\ \vdots & & \vdots \\ \cos\angle(X_n,X_0,X_1) & \hdots & \cos\angle (X_n,X_0,X_n) \end{array} \end{equation}

\noindent As a consequence of this and of Theorem \ref{theo.ec},  we can prove the following theorem answering Problem \ref{question.mengerplane}:

\begin{theorem}\label{theo.en} Let $(\mathcal X,\mathcal B,\angle)$ be an angle space. It can be conformally embedded in $\mathbb E^2$ if and only if every tetragon $\{X_0,X_1,X_2,X_3\}\subset\mathcal X$ admits an Euclidean-compatible distance function and  satisfies:

\noindent \emph{(Adjacent Angles Property)} Among the three angles $\angle(X_1,X_0,X_2)$, $\angle(X_1,X_0,X_3)$, $\angle(X_2,X_0,X_3)$ one of them must be the sum of the other two or $2\pi$ minus this sum (we have defined angles to be in the interval $[0,\pi]$).

% one of $\angle(B,A,C)\pm\angle(A,B,C)$ must equal $\angle(A,C,B)$ or  $2\pi-\angle(A,C,B)$ (compare with (EC2)).

\end{theorem}

\begin{proof} $(\mathcal X,\mathcal B,\angle)$ can be endowed with a unique Euclidean-compatible distance $d$, up to multiplication by a constant. It can be conformally embedded in $\mathbb E^2$ if and only if $(\mathcal X,d)$ can be isometrically embedded in $\mathbb E^2$. We can see that forall  $\{X_0,X_1,X_2\}\subset \mathcal X$, the determinant appearing in \eqref{eq.cosenoso} for $n=2$ is $1-\cos^2\angle(X_1,X_0,X_2)$ which is obviously greater or equal to $ 0$. On the other hand the same determinant for $n=3$ equals 0 if and only if  
$$\cos^2\angle(X_1,X_0,X_2)+\cos^2\angle(X_1,X_0,X_3)+\cos^2\angle(X_2,X_0,X_3)-$$
$$-2\cos\angle(X_1,X_0,X_2)\cdot \cos\angle(X_1,X_0,X_3)\cdot \cos \angle(X_2,X_0,X_3)-1=0 $$

\noindent which is equivalent to 
$$\big[\cos(\angle (X_1,X_0,X_2)+\angle (X_2,X_0,X_3))-\cos\angle(X_1,X_0,X_3)\big]\cdot$$
$$\cdot \big[\cos(\angle(X_1,X_0,X_2)-\angle(X_2,X_0,X_3))-\cos\angle(X_1,X_0,X_3)\big]=0$$
\noindent and, in turn, to the Adjacent Angles Property, since we must be in  one of the following cases
$$\begin{array}{l}\angle(X_1,X_0,X_2)+\angle(X_2,X_0,X_3)=\angle(X_1,X_2,X_3)\\
\angle(X_1,X_0,X_2)+\angle(X_2,X_0,X_3)=2\pi-\angle(X_1,X_2,X_3)\\
| \angle(X_1,X_0,X_2)-\angle(X_2,X_0,X_3)|=\angle(X_1,X_2,X_3)\end{array}$$
\end{proof}

Note that, together, the \emph{Axiom of Collinearity} in Definition \ref{defi.angle} and Condition (EC2) in Theorem \ref{theo.ec} form a stronger version of the \emph{Adjacent Angle Property} for the restricted cases in which  three of the points in $\{X_0,X_1,X_2,X_3\}$ are imposed to be collinear.

%%%%%%%%%%%%%%%%%%%%%%%%%%%%%%%%%%%%%%%%%%%%%%%%%%%%%%%%%%
\section{Final Comments} \label{sect.final}

%The intention of Theorem  \ref{theo.ec} is to provide a list of conditions that an angle space should satisfy to admit a Euclidean-compatible distance function.  Let us think about a space $\mathcal X$ consisting in 6 points (hexagon), in such a way that no three of them are collinear. This space has 60 angles. It also has $\binom{6}{3}$ trigons, $\binom{6}{4}$ tetragons and $\binom{6}{5}$ pentagons, each of them with the corresponding conditions required in (EC1), (EC2), (EC3) and (EC4). It is now clear that some of these conditions must be redundant. An improved statement of this theorem would be desirable, not only for aestetic reasons, but for a better understanding of angle spaces.

Given a metric space $(\mathcal X,d)$ equation \eqref{eq.ec} gives a unique candidate $\angle$ for a Euclidean-compatible angle function. But this function may fail to satisfy the Axiom of Collinearity.  It can be proved straightforwardly that this $\angle$ is an angle function if and only if for every tetragon $\{A,B,C,D\}\subset \mathcal X$ such that $D$ is between $B,D$ and $A, B$ and $C$ are not collinear  (with the relation of betweenness induced by $d$), the relation  of Stewart's Theorem holds, that is:
$$ d(A,C)^2\cdot d(B,D)+d(A,B)^2\cdot d(C,D) =d(B,C)\cdot (d(A,D)^2+d(B,D)\cdot d(C,D))$$

In general,  the characterizations given in Theorems \ref{theo.ec} and \ref{theo.en} can only be verified, in practice, for spaces $\mathcal X$ with a finite and small number of points. So it would be desirable to find alternative characterizations of angle spaces conformally embeddable in $\mathbb E^2$ (and $\mathbb E^n$) more applicable for other families of angles spaces of interest.

One of those familes of interest could be, for instance, the class $\mathcal T$ of angle spaces $(\mathcal X,\mathcal B, \angle)$ such that  $\mathcal X$ is additionally a topological space homeomorphic to the unit circle $S^1$. There is a correspondence between \emph{plane Jordan curves} and elements in $\mathcal T$ and between \emph{knots} and elements in $\mathcal T$. Understanding the angle space structure of these spaces may be very important. Note, for example, that a well known old and unsolved conjecture, the \emph{Inscribed Square Problem} (see \cite{M.IS, T.IS}), can be stated in the following terms:

\begin{conj} Let  $(\mathcal X,\mathcal B, \angle)$ be an angle space such that $\mathcal X$ is a topological space which is homeomorphic and conformally equivalent to a plane Jordan curve. Then there exists a tetragon $\{A,B,C,D\}$ such that  $\angle(A,B,C)=\angle(B,C,D)=\angle(C,D,A)=\angle(D,A,B)=\pi/2$ and $\angle(B,A,C)=\angle(C,A,D)=\pi/4$.

\end{conj}

 Problems of isometric embeddability of finite spaces  in Euclidean and hyperbolic spaces have recently received a lot of attention due to their importance in machine learing (see \cite{MAT}).  The concept of angle function may contribute to this class of problems. Theorem \ref{theo.en} can be easily generalized to characterize those angle spaces that can be conformally embedded in $\mathbb E^n$.  It is also possible to consider problems of conformal embeddability in other ambient spaces with a notion of angle, such as hyperbolic spaces or other manifolds.

In view of this and to finish this article, we  propose three problems for the future:

%\begin{problem} Decide whether or not conditions (EC1), (EC2), (EC3) and (EC4) are redundant or can be simplified.\end{problem}

\begin{problem} Let $\mathcal X$ be a topological space homeomorphic to $S^1$. Find necessary and sufficient conditions for the angle space $(\mathcal X,\mathcal B,\angle)$ to be conformally equivalent to a plane Jordan curve, easier to verify than the ones in Theorem \ref{theo.en}.

\end{problem}

\begin{problem} Characterize those angle spaces  that can be conformally embedded in a sphere in $\mathbb E^3$.

\end{problem}

\begin{problem} Characterize those angle spaces that can be conformally embedded in the hyperbolic plane.

\end{problem}

\section*{Acknowledgements}

The author want to thank the referees for their helpful suggestions and comments, that improved remarkably the first version of this paper.

% ------------------------------------------------------------------------
\end{document}